\documentclass{aptpub}

\usepackage{color}
\usepackage{amssymb}
\usepackage{amsfonts}
\usepackage{amsmath}
\usepackage{mathrsfs}

\authornames{Damien Bankowski, Claudia Kl\"uppelberg and Ross Maller}
\shorttitle{Ruin probability of the genOU process}

\numberwithin{equation}{section}

\newcommand{\bexam}{\begin{example}\rm}
\newcommand{\eexam}{\end{example}}
\newcommand{\bproof}{\begin{proof}}
\newcommand{\eproof}{\end{proof}}

\newcommand{\beao}{\begin{eqnarray*}}
\newcommand{\eeao}{\end{eqnarray*}\noindent}

\newcommand{\beam}{\begin{eqnarray}}
\newcommand{\eeam}{\end{eqnarray}\noindent}

\newcommand{\barr}{\begin{array}}
\newcommand{\earr}{\end{array}}

\def\bbr{{\Bbb R}}
\def\bbn{{\Bbb N}}

\newcommand{\barZ}{{\overline Z}}
\newcommand{\whZ}{{\widehat Z}}
\newcommand{\whM}{{\widehat M}}
\newcommand{\whQ}{{\widehat Q}}
\newcommand{\whL}{{\widehat L}}
\newcommand{\whX}{{\widehat X}}
\newcommand{\whc}{{\widehat c}}

\newcommand{\barL}{{\overline L}}

\newcommand{\veps}{\varepsilon}

\newcommand{\al}{{\alpha}}

\newcommand{\la}{{\lambda}}

\newcommand{\si}{{\sigma}}

\newcommand{\wh}{\widehat}

\newcommand{\halmos}{\quad\hfill\mbox{$\Box$}}

\newcommand{\cadlag} {c\`{a}dl\`{a}g\ }
\newcommand{\Levy} {L\'{e}vy~}

\newcommand{\ud}{\mathrm{d}}

\newcommand{\R}{\mathbb{R}}
\newcommand{\rmd}{{\rm d}}

\begin{document}

\title{On the Ruin Probability of the Generalised \\ Ornstein-Uhlenbeck Process in the Cram\'er Case$^\dagger$\footnote{$^\dagger${This
research was partially supported by ARC grant DP1092502.}}}

\authorone[Australian National University]{Damien Bankovsky}
\addressone{Mathematical Sciences Institute, Australian National University, Canberra, Australia, email: Damien.Bankovsky@anu.edu.au }
\authortwo[Technische Universit\"at M\"unchen]{Claudia Kl\"uppelberg}
\addresstwo{Center for Mathematical Sciences, and Institute for Advanced Study, Technische Universit\"at M\"unchen,  85747 Garching, Germany,
email: cklu@ma.tum.de}
\authorthree[ Australian National University]{Ross Maller}
\addressthree{Mathematical Sciences Institute, and
School of Finance and Applied Statistics, Australian National University, Canberra, Australia, email: Ross.Maller@anu.edu.au}

\begin{abstract}
For a bivariate \Levy process $(\xi_t,\eta_t)_{t\ge 0}$ and initial value $V_0$ define the Generalised Ornstein-Uhlenbeck (GOU) process
\[
V_t:=e^{\xi_t}\Big(V_0+\int_0^t e^{-\xi_{s-}}\ud \eta_s\Big),\quad t\ge0,\]
and the associated stochastic integral process
\[Z_t:=\int_0^t e^{-\xi_{s-}}\ud \eta_s,\quad t\ge0.\]
Let
$T_z:=\inf\{t>0:V_t<0\mid V_0=z\}$ and
$\psi(z):=P(T_z<\infty)$ for $z\ge 0$ be
the ruin time and infinite horizon ruin probability of the GOU.
Our results extend previous work of
Nyrhinen (2001) and others
to give  asymptotic estimates for $\psi(z)$ and the distribution
of $T_z$ as $z\to\infty$,
under very general, easily checkable, assumptions,
when $\xi$ satisfies a Cram\'er condition.
\end{abstract}



\section{Introduction}\label{s1}

Let $(\xi,\eta)=(\xi_t,\eta_t)_{t\ge 0}$ be a bivariate \Levy process on a filtered complete probability space $(\Omega, \mathscr{F},\mathbb{F},P)$ and define
a generalised Ornstein-Uhlenbeck (GOU) process by
\begin{equation}\label{GOU definition}
V_t:=e^{\xi_t}\Big(V_0+\int_0^t e^{-\xi_{s-}}\ud \eta_s\Big),\ t\ge 0,
\end{equation}
and the associated stochastic integral process $Z=(Z_t)_{t\ge 0}$ by
\begin{equation}\label{Z definition}
Z_t:=\int_0^t e^{-\xi_{s-}}\ud \eta_s.
\end{equation}
$V_0$ is a random variable (r.v.), not necessarily independent of  $(V_t)_{t>0}$.
To avoid trivialities, assume that neither $\xi$ nor $\eta$ are
identically zero.

Such processes have attracted attention over the last decade as continuous time analogues of solutions to stochastic recurrence equations (SRE); cf. Carmona, Petit and Yor~\cite{CarmonaPetitYor97,CarmonaPetitYor01},
Erickson and Maller~\cite{EricksonMaller05}.
The link between SREs and the GOU was made in
de Haan and Karandikar~\cite{deHaanKarandikar89}.
GOU processes turn up naturally in
 stochastic volatility models (e.g., the continuous time GARCH model of Kl\"uppelberg, Lindner and Maller~\cite{klm:2004}), but most prominently as insurance risk models
for perpetuities in life insurance or
when the insurance company receives some stochastic return on investment; such investigations started with Dufresne~\cite{Dufresne90} and Paulsen~\cite{Paulsen93}.
More references are given later.

This paper is intended to fill a gap left
between Bankovsky~\cite{Bankovsky09} and
Bankovsky and Sly~\cite{BankovskySly08}, where more details on
the insurance background can be found.
Define
$$
T_z:=\inf\{t>0:V_t<0 \mid V_0=z\}, \ z\ge 0,
$$
(with the convention throughout that  $\inf\emptyset=\infty$),
and let
\begin{equation}\label{definition: ruin prob}
\psi(z):=P\Big(\inf_{t>0}V_t<0\mid V_0=z\Big)
=P\Big(\inf_{t>0}Z_t<-z\Big)=P\left(T_z<\infty\right),\ z\ge0,
\end{equation}
be the \emph{infinite horizon ruin probability} for the GOU.
Note that $\psi(z)$ is a nonincreasing function of $z$,
and we can ask how fast it decreases as $z\to\infty$.

Our main result, Theorem~\ref{theorem: my Nyrhinen result},
provides a very general asymptotic result for $\psi(z)$ as $z\to\infty$
for the case when
$\lim_{t\to\infty}Z_t$ exists as an a.s. finite r.v.
and shows that, under a Cram\'er-like condition on $\xi$,
$\psi(z)$ decreases approximately like a power law.
This is an extension of a similar asymptotic result of
Nyrhinen \cite{Nyrhinen01}, who, like us,
utilises a discrete time result  of
Goldie \cite{goldie:1991} for proof.
We use more recent developments in the theory of discrete time
perpetuities and the continuous time GOU to update  Nyrhinen's results.
In Section \ref{Examples} we provide some examples
which cannot be dealt with by the prior results
but satisfy the conditions of our theorem.

To conclude this introduction,
we describe some previous literature relating to the GOU and its ruin
probability, beginning with those papers which examine the GOU in
its full generality. The process appears implicitly in the work of
de Haan and Karandikar~\cite{deHaanKarandikar89} as a continuous
generalisation of an SRE. Basic properties
are given by Carmona et al. \cite{CarmonaPetitYor01}. A general
survey of the GOU and its applications is given by Maller, M\"uller
and Szimayer \cite{MallerMullerSzimayer07}.
Exact conditions for no ruin
($\psi(z)=0$ for some $z\ge 0$) are given by
Bankovsky and Sly \cite{BankovskySly08} whilst conditions for
certain ruin
($\psi(z)=1$  for some $z\ge 0$)
are examined by Bankovsky \cite{Bankovsky09}.

The study of the GOU is closely related to the study of integrals of
the form $Z,$ defined in (\ref{Z definition}).
It is shown in Lindner and Maller \cite{lindner:maller:2005} that
stationarity of $V$ is related to convergence of a
stochastic integral constructed from $(\xi,\eta)$ in a similar way to $Z.$

Among the few papers dealing with $Z$ in its full generality,
Erickson and Maller \cite{EricksonMaller05} give necessary and
sufficient conditions for the almost sure convergence of $Z_t$ to a
r.v. $Z_\infty$ as $t\rightarrow\infty$,
and Bertoin, Lindner and Maller
\cite{BertoinLindnerMaller07} present necessary and sufficient
conditions for the continuity of the distribution of
$Z_\infty$, when  it exists.
Fasen~\cite{fasen:2009}, using point process methods, gives an account
of the extremal behaviour of a GOU process.

There are a larger number of papers dealing with $V$ and $Z$ when
$(\xi,\eta)$ is subject to restrictions. We discuss a selection of those
papers which are relevant to ruin probability.
Harrison \cite{Harrison77} presents results on the ruin probability
of $V$ when $\xi$ is a linear deterministic function and $\eta$ is a
\Levy process with finite variance.
His approach is based on an exponential martingale argument,
which corresponds to the Cram\'er case.
The heavy-tailed case is investigated in
Kl\"uppelberg and Stadtm\"uller~\cite{Klu:Stadtmueller:98} and
extended by Asmussen~\cite{Asmussen:1998}.
See also Maulik and Zwart~\cite{MaulikZwart} and Konstantinides and Mikosch \cite{KonstantinidesMikosch05}.

Paulsen \cite{Paulsen93}
generalises Harrison's results, and presents new ruin probability
results for $V$, when $\xi$ and $\eta$ are independent with
finite activities.
This independent case is also treated in Kalashnikov and Norberg
\cite{KalashnikovNorberg02} and Paulsen \cite{Paulsen98,Paulsen02}.
Chiu and Yin \cite{ChiuYin04} generalise some of Paulsen's results
to the case in which $\eta$ is a jump-diffusion process.
Cai \cite{Cai04} and Yuen et al. \cite{NgYuenWang04}
present results when $\eta$ is a compound Poisson process.

Most relevant works containing restrictions on $(\xi,\eta)$ focus on the case
when $Z_t$ converges to $Z_\infty$ as $t\rightarrow\infty$; cf.
 Yor \cite{Yor01} and Carmona et al.
\cite{CarmonaPetitYor97}. Gjessing and Paulsen
\cite{GjessingPaulsen97} study the distribution of $Z_\infty$ when
$\xi$ and $\eta$ are independent with finite activity,
and obtain exact distributions in some special cases. Hove and
Paulsen \cite{PaulsenHove99} use Markov chain Monte Carlo methods
to find the distribution of $Z_\infty$ in some special cases. Kl\"uppelberg and Kostadinova
\cite{KluppelbergKostadinova08} and Brokate et al.
\cite{Malleretal08} provide results on the tail of the distribution of
$Z_\infty$ when
$\eta$ is a compound Poisson process plus drift, independent of $\xi$.

\section{Main Results}\label{s2}

Our main results apply under a
 Cram\'er-like condition on $\xi$: assume that
\begin{equation}\label{cramer}
Ee^{-w\xi_1}=1 \mbox{ for some }  w>0.
 \end{equation}
The following consequences of  \eqref{cramer}
are well known and easily verified.
Condition \eqref{cramer} implies that $E\xi_1$ is well defined,
with $E\xi_1^-<\infty$, $E\xi_1^+\in(0,\infty]$,
and $E\xi_1\in (0, \infty]$, and  so
$\lim_{t \to \infty} \xi_t=\infty$ a.s.
Further,
$Ee^{-\alpha \xi_1}$ is finite and nonzero for all $\alpha\in[0,w]$,
and
$c(\alpha):= \ln Ee^{-\alpha \xi_1}$
is finite at least for all $\alpha\in[0,w)$.
The derivatives $c'(\alpha)$ and  $c''(\alpha)$ are
finite at least for all $\alpha\in[0,w)$,
and $c''(\alpha)\in(0,\infty]$ for all $\alpha\ge 0$.
So $c(\alpha)$ is strictly convex for $\alpha\in[0,\infty)$
and $\mu^*:= c'(w)=-E[\xi_1e^{-w\xi_1}]\in(0,\infty]$.

We will need the {\em Fenchel-Legendre transform} of $c$, defined as
\beam\label{fenchel}
c^*(v):=\sup\{\alpha
v-c(\alpha):\alpha\in\mathbb{R}\}, \ v\in\mathbb{R}.
\eeam
Next, let
\begin{equation}\label{t0}
\alpha_0:= \sup\left\{\alpha\in\mathbb{R}:c(\alpha)<\infty,
E|Z_1|^\alpha<\infty
\right\}\in[0,\infty],
\end{equation}
and define the constant
\begin{equation}\label{xdef}
x_0:=\lim_{\alpha\rightarrow
 \alpha_0-} (1/c'(\alpha))\in[0,\infty].
 \end{equation}
A distribution is \emph{spread out}
if it has a convolution power with an absolutely continuous component.

\begin{thm}\label{theorem: my Nyrhinen result}
Suppose that the following conditions hold:\\
{\rm Condition A: } $\psi(z)>0$ for all $z\ge 0$,\\
{\rm Condition B: } there exists $w>0$ such that $Ee^{-w\xi_1} = 1$ (i.e. \eqref{cramer} holds),\\
{\rm Condition C: }
there exist $\varepsilon>0$
and $p,q>1$ with $1/p+1/q=1$ such that
\begin{equation}\label{equation: first moment condition}
E[e^{-\max\{1,w+\varepsilon\}p\xi_1}]<\infty\quad\mbox{and}\quad
E[|\eta_1|^{\max\{1,w+\varepsilon\}q}]<\infty.
\end{equation}
Then $0\le x_0<1/\mu^*<\infty$, the function
\[
R(x):= \left\{ \begin{array}{ll}
xc^*(1/x) & \textrm{for~}x\in(x_0,1/\mu^*), \\
w &\textrm{for~}x\ge 1/\mu^*,
\end{array} \right.\,
\]
is finite and continuous on $(x_0,\infty)$ and strictly decreasing
on $(x_0, 1/\mu^*)$, and we have
\begin{equation}\label{equation: Nyrhinen result 1}
\lim_{z\rightarrow\infty}(\ln z)^{-1} \ln P(T_z\le x\ln z) = -R(x)
\end{equation}
for every $x>x_0.$
In addition,
\begin{equation}\label{equation: Nyrhinen result 2}
\lim_{z\rightarrow\infty}(\ln z)^{-1} \ln \psi(z)=-w.
\end{equation}
If, further, the distribution of $\xi_1$ is spread out, then there
exist constants $C_->0$ and $\kappa>0$ such that
\begin{equation}\label{equation: Goldie
result}z^w\psi(z)=C_- + o(z^{-\kappa})~~\mathrm{as~}z\rightarrow\infty.
\end{equation}
\end{thm}

\begin{rem}\rm
(i)  $\psi(z)>0$ for all $z\ge 0$ is of course a logical assumption
to make in the context of Theorem~\ref{theorem: my Nyrhinen result},
though not necessarily easy to verify.
Necessary and sufficient conditions for it
in terms of the L\'evy measure of $(\xi,\eta)$
are given  in \cite{BankovskySly08}.
The moment conditions in Theorem~\ref{theorem: my Nyrhinen result}
are also easily expressed in terms of the L\'evy measure of
$(\xi,\eta)$, cf. Sato \cite{sato:1999}, p.~159.
They imply that
$E[\sup_{0\le t\le 1}\left|Z_t\right|^{\max\{1,w+\veps \}}]<\infty$
(see Lemma \ref{sup proposition} below).
We also have $E[\ln(\max\{1,|\eta_1|\}]<\infty$ in  Theorem~\ref{theorem: my Nyrhinen result},
and $\lim_{t \to \infty} \xi_t=\infty$ a.s., so
$Z_t$ converges a.s. to a finite r.v.
$Z_\infty$ as $t\rightarrow\infty$
by Proposition 2.4 of \cite{lindner:maller:2005}
or Theorem 2 of \cite{EricksonMaller05}.\\
(ii)
Let
$\barZ_t:= Z_t-\inf_{0\le s\le t}Z_s$ be the process
reflected in its minimum, and set
\begin{equation}\label{MQL}
(M,Q,\barL):=\big(e^{-\xi_1},Z_1,-e^{\xi_1}\barZ_1\big).
 \end{equation}
Then the value $C_-$ in
(\ref{equation: Goldie result}) is given by the formula
in (2.19) of Goldie  \cite{goldie:1991}, namely
\begin{equation}C_-=\frac{1}{w
\mu^*}E\Big[
\Big(Q+M\min\big\{\barL,\inf_{t>0}Z_t\big\}^-\Big)^w-
\Big(\big(M\inf_{t>0}Z_t\big)^-\Big)^w\Big].
\end{equation}
When $\xi$ and $\eta$
are independent, it was pointed out by Paulsen \cite{Paulsen02}
that this constant can be written in a slightly different form,
which, by
Theorem 4 of \cite{BankovskySly08}, is also true
in the dependent case. Namely, let
$G(z):=P(Z_\infty\le z)$,
$h(z):=E[G(-V_{T_z})\mid T_z<\infty]\in[0,1]$,
and
$h:=\lim_{z\rightarrow\infty} h(z)$.
Then
\[
C_-=\frac{1}{w\mu^* h}
E\Big[\left(\left(Q+MZ_\infty\right)^-\right)^w-
\left(\left(MZ_\infty\right)^-\right)^w\Big].
\]
(iii)  The requirement that $\xi_1$ is
spread out can be replaced
with the less restrictive requirement that $\xi_T$ be spread out, where
$T$ is uniformly distributed on $[0,1]$ and independent of $\xi.$
We omit details of this, which can be carried out as in
\cite{Paulsen02}.
\end{rem}

\section{Examples}\label{Examples}

In this section we provide examples of \Levy processes for
which Conditions A, B and C of Theorem~\ref{theorem: my Nyrhinen result}
are satisfied. 
Note that conditions B and C only involve the marginal processes $\xi$ and $\eta$ and they apply to all examples treated in the literature so far; cf. Kl\"uppelberg and Kostadinova~\cite{KluppelbergKostadinova08} for detailed references.
The only condition which may involve dependence between $\xi$ and $\eta$ is
Condition A.

We denote the characteristic triplet of $(\xi,\eta)$ by
$((\tilde{\gamma}_\xi,\tilde{\gamma}_\eta),\Sigma_{\xi,\eta},\Pi_{\xi,\eta}).$
The characteristic triplet of the marginal process $\xi$ is denoted by  $(\gamma_\xi,\sigma_\xi^2,\Pi_\xi),$ where
\begin{equation}\label{first 2 dim to 1 dim equation}
\gamma_\xi=\tilde\gamma_\xi+
\int_{\{|x|< 1\}\cap\{x^2+y^2\ge 1\}}x\Pi_{\xi,\eta}(\ud(x,y)),
\end{equation}
and $\sigma_\xi^2$ is the upper left entry in the matrix $\Sigma_{\xi,\eta}$.
Similarly for $\eta$.
The random jump measure and Brownian motion
components of $(\xi,\eta)$ will be denoted respectively by
$N_{\xi,\eta}$ and $(B_\xi,B_\eta)$;
see Section 1.1 of \cite{BankovskySly08} for further details.

\bexam [Bivariate compound Poisson process with drift]\label{3.1}
\\
Let $(N_t)_{t\ge0}$ be a Poisson process with intensity $\la>0$,
and, independent of it, $(X_i,Y_i)_{i\in\bbn}$ an iid sequence of random 2-vectors.
For $\gamma_\xi,\gamma_\eta\in\bbr$ set
\[
(\xi_t,\eta_t):=(\gamma_\xi,\gamma_\eta)\,t+\sum_{i=1}^{N_t}(X_i,Y_i),\quad t\ge0,
\]
with
$E|X_1|<\infty$ and $\lambda$, $\gamma_\xi$ and $EX_1$ such that $\gamma_\xi+\lambda EX_1>0$.
For this process,
\[
c(\al)=\ln Ee^{-\alpha\xi_1}=-\al \gamma_\xi-\lambda \big(1-Ee^{-\al X_1}\big) < \infty
\]
for $\al\in\mathbb{R}$ such that $Ee^{-\alpha X_1}$ is finite,
with $c'(0)=-\gamma_\xi-\lambda EX_1<0$.

We consider the special case where  $(X_1,Y_1)$ is bivariate Gaussian with mean
$(m_X,m_Y)$ and positive definite covariance matrix
\[ \Sigma_{X,Y}:=
\left(
 \begin{array}{cc}
\sigma_X^2 & \sigma_{X,Y} \\
\sigma_{X,Y} & \sigma_Y^2  \\
\end{array}
\right).
\]
Then Condition C obviously  holds.
For Condition B, note that
\begin{equation}\label{ca}
c(\alpha)=
-\al \gamma_\xi - \la \big(1-e^{-m_X \al+\si_X^2\al^2/2}\big) \to \infty\,\mbox{ as } \alpha\to\infty.
\end{equation}
Consequently, a Lundberg coefficient exists and Condition B is satisfied.
To establish Condition A we note that $(\xi,\eta)$ is a finite variation process and invoke Remark 2(2) of \cite{BankovskySly08},
also using the notation from that paper.
In fact, by that Remark 2(2), $\psi(z)=0$ for some $z>0$ would imply that
$P_{X,Y}(A_3)=P(X_1\le 0, Y_1\le 0)=0$, which obviously
is not the case. So Condition A holds.
\eexam

\bexam A Brownian motion with drift, i.e., with
\[
(\xi_t,\eta_t)=(\gamma_\xi,\gamma_\eta)\,t+(B_{\xi,t},B_{\eta,t}),\quad t\ge0,
\]
where $\gamma_\xi>0$ and
$(B_\xi,B_\eta)_t$ is bivariate Brownian motion with mean 0 and positive definite
covariance matrix, is easily seen to satisfy Conditions A, B, C.
\eexam

\bexam [Jump diffusion $\xi$ and Brownian motion $\eta$]\\
Let $(B_t)_{t\ge0}$ be Brownian motion with mean zero and variance $\si^2$, $(N_t)_{t\ge0}$ a
Poisson process with intensity $\la>0$, and $(X_i)_{i\in\bbn}$ iid r.v.s, all independent.
Set
$$
(\xi_t,\eta_t) = (\gamma_\xi,\gamma_\eta)t + \big(B_t+\sum_{i=1}^{N_t} X_i, B_t\big),\quad t\ge0,
$$
where $\gamma_\xi>0$, and
assume that $\gamma_\xi+\lambda EX_1>0$.
Condition A holds, since the Gaussian covariance matrix of $(\xi,\eta)$ is of the form
\beam\label{matrix}
 \Sigma_{\xi,\eta}:=
\left(
 \begin{array}{cc}
\sigma^2+ \lambda EX_1^2 & 1 \\
1 & \sigma^2  \\
\end{array}
\right),
\eeam
and, hence, is not of the form excluded by Theorem~1 of \cite{BankovskySly08}.
Moreover, $c(\al)$ is the same as in \eqref{ca} with the addition of a term
$ \al^2\sigma^2/2$,
so again $c'(0)= - \gamma_\xi - \la EX_1 <0$.
\\[2mm]
(a) \, Now assume that $X_1$ is, as in the
Merton model, normally distributed with mean $m_X$ and variance $\si_X$.
Then
Conditions B and C are  satisfied just as in Example \ref{3.1}.
\\[2mm]
(b) \,
The picture changes slightly when we consider Laplace distributed $X$
with density $f(x)=\rho e^{-\rho |x|}/2$ for $x\in\bbr$, $\rho>0$.
Then $Ee^{-\al X}=\rho\left((\rho+\al)^{-1}+(\rho-\al)^{-1}\right)/2$
for $-\rho<\al<\rho$ with  singularities at $-\rho$ and $\rho$.
Moreover,
$$
c'(\al)=-\gamma_\xi +\al\sigma^2 +\la\frac{\rho}2\Big(\frac1{(\rho-\al)^2}
-\frac1{(\rho+\al)^2}\Big),
$$
implying that $c'(0)=-\gamma_\xi<0$.
So a Lundberg coefficient $w>0$ exists.
Since the normal r.v. $B_1$ has absolute moments of every order,
for Condition C to hold it suffices that $w<\rho$, which is guaranteed,
since $\rho$ is a singularity of $c$.
\eexam

\bexam [Subordinated Brownian motion $\xi$ and spectrally positive $\eta$]\\
Let $(B_t)_{t\ge 0}$ be a standard Brownian motion and $(S_t)_{t\ge 0}$
a driftless subordinator with $\Pi_S\{\R\}=\infty$.
For constants $\mu$, $\gamma_\xi$, $\gamma_\eta$,
define
$$(\xi_t,\eta_t) = (\gamma_\xi,\gamma_\eta) t + (B(S_t)+\mu S_t,S_t),\quad t\ge0.$$
Subordinated Brownian motions play an important role in financial modeling; cf. Cont and Tankov~\cite{CT}, Ch. 4.
The bivariate process above has joint Laplace transform
\beao
e^{(\alpha_1\gamma_\xi+\alpha_2\gamma_\eta)t}
E[e^{\alpha_1 (B(S_t)+\mu S_t)+\alpha_2S_t} ]
&=&
e^{(\alpha_1\gamma_\xi+\alpha_2\gamma_\eta)t}
 E[ e^{\Psi_B(\alpha_1) S_t + (\alpha_1\mu +\alpha_2 )S_t}]\\
&=&
e^{t[\Psi_S(\Psi_B(\alpha_1)+\alpha_1\mu+\alpha_2)+
\alpha_1\gamma_\xi+\alpha_2\gamma_\eta]},
\eeao
where $\Psi_B$ and $\Psi_S$ are the Laplace exponents of $B$ and
$S$, respectively.  Thus $\Psi_B(\alpha)=-\alpha^2/2$.
By setting $\alpha_2=0$ and $t=1$ we obtain
\beao
c(\alpha)= \ln E e^{-\al\xi_1}
= \Psi_S(\Psi_B(-\alpha) -\alpha\mu  )-\alpha \gamma_\xi
= \Psi_S(-\al^2/2 -\alpha\mu)-\alpha \gamma_\xi.
\eeao
Consider the variance gamma model with parameters $c,\la>0$, where $S$ is a gamma
subordinator with L\'evy density
$\rho(x)=cx^{-1}e^{-\la x}$ for $x>0$ and
Laplace transform
$Ee^{-uS_t} = (1+u/\la)^{-ct}$.
Assume $\gamma_\xi+c\mu/\lambda>0$ and $\gamma_\eta\le 0$.
Now, $\Psi_S(u)=-c\ln(1-u/\la)$,
giving
$$
c(\al)=-\alpha \gamma_\xi-c\ln\Big(1+\frac{\al\mu}{\la}-\frac{\al^2}{2\la}\Big).
$$
$c(\al)$ is well defined for $\alpha\in(\mu-\sqrt{\mu^2+2\la},\mu+\sqrt{\mu^2+2\la})$, which includes 0, and $c'(0)=-\gamma_\xi-c\mu/\la<0$.
Then, since $c(\mu+\sqrt{\mu^2+2\la})=\infty$, the Lundberg coefficient $w$ exists.

In order to check Condition A, we have, in the notation
of Theorem 1 of \cite{BankovskySly08},
$\Pi_{\xi,\eta}(A_2)=\Pi_{\xi,\eta}(A_3)=0$,
since $\eta$ has only positive jumps, and
$\theta_2=0$.
Now with $u\ge 0$,
$A_4^u=\{x\le 0,y\ge 0:y<u(e^{-x}-1)\}=\{x\ge 0,y\ge 0:y<u(e^x-1)\}$.
Since $\Pi_\eta(\R)=\infty$, $\eta$ has jumps arbitrarily close to 0, and we have
$\Pi_{\xi,\eta}(A_4^u)>0$ for $u>0$, while  $\Pi_{\xi,\eta}(A_4^0)=0$.
Thus $\theta_4:= \inf\{u\ge 0: \Pi_{\xi,\eta}(A_4^u)>0\}=0$.
There is no Gaussian component, so $\sigma^2_\xi=0$, which puts us in the situation
of the second item of Theorem 1 of \cite{BankovskySly08},
and to verify that $\psi(z)>0$ for all $z\ge 0$ we only need
(since $\theta_2=\theta_4=0$)
\begin{equation}\label{gcr}
g(0)=\widetilde \gamma_\eta-\int_{x^2+y^2\le 1}y\Pi_{\xi,\eta}(\rmd x,\rmd y)<0.
\end{equation}
But by \eqref{first 2 dim to 1 dim equation},
\[
\widetilde \gamma_\eta=\gamma_\eta
-\int_{0\le y\le1,x^2+y^2>1}y\Pi_{\xi,\eta}(\rmd x,\rmd y)\le\gamma_\eta,
\]
thus
$ g(0)<\gamma_\eta\le 0$,
since we chose $\gamma_\eta\le 0$.
Hence  Condition A holds in this model.
\eexam

\section{Discrete Time Background and Preliminaries}\label{s4}

Our continuous time asymptotic results
will be transferred across from discrete time versions, and
our first task  in the present section is to show how $(V_t)_{t\ge0}$ can be
expressed  as a solution of one of two SREs, and give the associated discrete
stochastic series for $(Z_t)_{t\ge0}$.
Earlier papers  in this area also adopted this approach and we
will tap into some of their results in proving
Theorem~\ref{theorem: my Nyrhinen result}.

We begin by describing the discrete time setup we use. For
$n\in\bbn$ consider the SRE
\begin{equation}\label{definition: first difference equation for V}
Y_n=A_n Y_{n-1}+B_n,
\end{equation}
where $(A_n,B_n)_{n\in\bbn}$ is an iid sequence of $\mathbb{R}^2$-valued random
vectors independent of an initial r.v. $Y_0.$ The
recursion in \eqref{definition: first difference equation for V} can
be solved in the form
\begin{equation}\label{definition: first discrete series for V}
Y_n=Y_0\prod_{j=1}^n A_j+\sum_{i=1}^n\prod_{j=i+1}^n A_jB_i
\end{equation}
(with $\prod_{j=n+1}^n=1$).
{}From (\ref{GOU definition}) we can write, for $n\in\bbn$
\begin{equation}\label{formula: expansion of V}V_n=e^{\xi_n-\xi_{n-1}}\Big(e^{\xi_{n-1}}\big(V_0+\int_0^{n-1}e^{-\xi_{s-}}\ud \eta_s\big)\Big)
+e^{\xi_n}\int_{(n-1)+}^n e^{-\xi_{s-}}\ud \eta_s.
\end{equation}
Thus, if we let $Y_0=V_0$ and define the  $\R^2$-valued random vectors
\begin{equation}\label{definition: first A and B}
(A_n,B_n):=\Big(e^{\xi_n-\xi_{n-1}},e^{\xi_n}\int_{(n-1)+}^n
e^{-\xi_{s-}}\ud \eta_s\Big),
\end{equation}
then
$V_n$ satisfies (\ref{definition: first difference equation
for V}).
An alternative formulation considers for $n\in\bbn$ the SRE
\begin{equation}\label{definition: second difference equation for V}
Y_n=C_n Y_{n-1}+C_n D_n\,,
\end{equation}
where $(C_n,D_n)_{n\in\bbn}$ is an iid
sequence independent of $Y_0.$
The solution is
\begin{equation}\label{definition: second discrete series for V}
Y_n=Y_0\prod_{j=1}^n C_i+\sum_{i=1}^n\prod_{j=i}^n C_jD_i.
\end{equation}
Using (\ref{formula: expansion of V}) it is clear that
$V_n$ is a solution of (\ref{definition: second difference equation for V})
if we let $V_0=Y_0$ and define
\begin{equation}\label{definition: second C and D}
(C_n,D_n):=\Big(e^{\xi_n-\xi_{n-1}},e^{\xi_{n-1}}\int_{(n-1)+}^n
e^{-\xi_{s-}}\ud \eta_s\Big).
\end{equation}
Then it is easily verified that
\begin{equation}\label{definition: second discrete series for Z}
Z_n=\sum_{i=1}^n\prod_{j=1}^{i-1} C_j^{-1}D_i
\end{equation}
(with $\prod_{j=1}^0=1$).
Note that even when $\xi$
and $\eta$ are independent, the r.v.s $A_n$ and $B_n$ may
be dependent, and similarly for  $C_n$ and $D_n$.
But we have

\begin{lem}\label{lemma: iid}
\label{Appendix} $(A_n,B_n)_{n\in\bbn}$ and  $(C_n,D_n)_{n\in\bbn}$  are
iid sequences.
\end{lem}

\bproof
We begin by proving that the sequence $(C_n,D_n)_{n\in\bbn}$ is iid.
Fix $n\in\mathbb{N}$ and define the new \Levy process
$
(\bar{\xi}_s,\bar{\eta}_s) :=(\xi_{n-1+s}-\xi_{n-1},\eta_{n-1+s}-\eta_{n-1})$ for $s\ge0.
$
Thus
$(\bar{\xi}_s,\bar{\eta}_s)_{s\ge 0}=_D(\xi_s,\eta_s)_{s\ge 0}$.
Note that we can bring the term $e^{\xi_{n-1}}$ through the integral
sign in \eqref{definition: second C and D} and write
$D_n=\int_{(n-1)+}^n e^{-(\xi_{s-}-\xi_{n-1})}\ud \eta_s.$ $(\xi,\eta)$ has independent increments, so
$(C_n,D_n)$ is independent of $(C_m,D_m)$ for every $n\neq m.$
Now
\begin{eqnarray*}
(C_n,D_n)&=&\Big(e^{\xi_n-\xi_{n-1}},\int_{(n-1)+}^n e^{-\left(\xi_{s-}-\xi_{n-1}\right)}\ud\eta_s\Big)\\
&=&\Big(e^{\bar{\xi_1}},\int_{0+}^1e^{-\bar{\xi}_{s-}}\ud\bar{\eta}_s\Big)
 \, =_D \, \Big(e^{\xi_1},\int_{0+}^1e^{-\xi_{s-}}\ud\eta_s\Big)=(C_1,D_1).
\end{eqnarray*}
Thus we have proved that $(C_n,D_n)_{n\in\bbn}$ is an iid sequence. This implies that $(C_n,C_nD_n)$ is also an iid sequence,
and then $(A_n,B_n)_{n\in\bbn}$ is also an iid sequence since
\[
(C_n,C_nD_n)=\Big(e^{\xi_n-\xi_{n-1}},e^{\xi_n}\int_{(n-1)+}^n
e^{-\xi_{s-}}\ud \eta_s\Big)=(A_n,B_n).
\]\vskip-0.8cm
\halmos
\eproof

In order to directly access particular results from previous papers,
when discretizing $V$ we will use the approach via the recursion
(\ref{definition: first difference equation for V}) and the
sequence (\ref{definition: first discrete series for V}), whereas
when discretizing $Z$ we will use the approach via the series
(\ref{definition: second discrete series for Z}).
There has been significant attention paid to
sequences of the form (\ref{definition: first
discrete series for V}) and (\ref{definition: second discrete series
for Z}), and they are linked via the fixed point of the same SRE, see
Vervaat \cite{vervaat:1979} and Goldie and Maller \cite{goldie:maller:2000}.

Next we describe two important papers relating to the GOU
and its ruin time.
In them, $\xi$ and $\eta$
are general \Levy processes, possibly dependent. The
relevant papers are Nyrhinen \cite{Nyrhinen01} and  Paulsen
\cite{Paulsen02}, which are very closely related to Theorem
\ref{theorem: my Nyrhinen result}.

Nyrhinen \cite{Nyrhinen01} contains asymptotic ruin probability
results for the GOU, in which $(\xi,\eta)$ is allowed to be an
arbitrary bivariate \Levy process.
He discretizes the stochastic integral process $Z$ and deduces asymptotic results in the continuous time setting from  similar discrete time results.
We describe Nyrhinen's results in some
detail, and then make some comments.

Let $(M_n,Q_n,L_n)_{n\in\bbn}$ be iid random vectors
with $P(M>0)=1$ and
$(M,Q,L)\equiv(M_1,Q_1,L_1)$.
Define the sequence
$(X_n)_{n\in\bbn}$ by
\begin{equation}\label{definition: Nyrhinens series}
X_n=\sum_{i=1}^n\prod_{j=1}^{i-1}M_jQ_i+\prod_{j=1}^nM_jL_n, \ {\rm
with}\ X_0=0.
\end{equation}
For $u>0$ define the passage time
$\tau_u^X:=\inf\{n\in\bbn:X_n>u\}$
and the function
$c_M(\alpha):=\ln EM^\alpha$.
Assume there is a $w^+>0$ such that $EM^{w^+}=1$.
Define
\begin{equation}\label{t0+}
\alpha_0^+:= \sup\left\{\alpha\in\mathbb{R}:c_M(\alpha)<\infty,
~E|Q|^\alpha<\infty,
~E(ML^+)^\alpha<\infty\right\}\in[0,\infty].
\end{equation}
Also let
\begin{equation}\label{Nyybar}
\bar{y}:=\sup\Big\{y\in\mathbb{R}:P\big(\sup_{n\in\mathbb{N}}
X_n>y\big)>0\Big\}\in(-\infty,\infty].
\end{equation}
Nyrhinen provides asymptotic results for $X_n$ under the following
\\[2mm]
{\bf Hypothesis H}: Suppose that $0<w^+<\alpha_0^+\le\infty$ and
$\bar{y}=\infty.$\\[2mm]
Under Hypothesis H,
and assuming that
$P(M>1)>0,$ the
following quantities are well-defined:
$\mu^+:=1/{c_M'(w)}\in(0,\infty)$ and
$x_0^+:=\lim_{t\rightarrow \alpha_0^+-} (1/{c_M'(t)})\in[0,\infty)$.
Let $c_M^*(v)$
 be the Fenchel-Legendre transform of $c_M$ as in \eqref{fenchel}.
Define the function
$R:(x_0^+,\infty)\rightarrow\mathbb{R}\cup\{\pm\infty\}$ by
\[
R(x):= \left\{ \begin{array}{ll}
xc_M^*(1/x) & \textrm{for~}x\in(x_0^+,1/\mu^+), \\
w & \textrm{for~}x\ge 1/\mu^+.
\end{array} \right.\,
\]
In our
situation, $R$ is finite and continuous on $(x_0^+,\infty)$ and
strictly decreasing on $(x_0^+,1/\mu^+)$.

\begin{prop}\label{theorem: Nyrhinens main theorem}{\rm [Nyrhinen's main discrete results, \cite{Nyrhinen01}, Theorems~2 and~3]}\newline
Assume Hypothesis H. Then the following hold.\\
(i) \, For every $x>x_0$,
\begin{equation}\label{4.15}
\lim_{u\rightarrow\infty}(\ln u)^{-1}\ln P(\tau_u^X\le x\ln u)=-R(x)
\end{equation}
and
\begin{equation}\label{equation: Nyrhinens finite time}
\lim_{u\rightarrow\infty}(\ln u)^{-1}\ln P(\tau_u^X<\infty)=-w.
\end{equation}
(ii) \, If the distribution of $\ln M$ is spread out,
there are constants  $C_+>0$ and $\kappa>0$ such that
\begin{equation}\label{equation: Nyrhinens Goldie}
u^{w^+}P(\tau_u^X<\infty)=C_+ +o(u^{-\kappa}), \ {\rm as}\ u\to\infty.
\end{equation}
\end{prop}

$C_+$ can be obtained from the formula in Theorem 6.2 and (2.18) of Goldie
\cite{goldie:1991}. Nyrhinen continues in his Theorem 3 to give
equivalences for the condition $\bar{y}=\infty$, but they are difficult to verify, as he admits.  We discuss these more fully later.

Nyrhinen's continuous result is obtained by applying his discrete
results to the  case
 \begin{eqnarray}\label{definition: M allocation}
(M_n, Q_n)&=&
\Big(e^{-(\xi_n-\xi_{n-1})},
e^{\xi_{n-1}}\int_{(n-1)+}^ne^{-\xi_{s-}} \ud \eta_s\Big)
=(C_n^{-1}, D_n)
\quad  {\rm (cf.\ \eqref{definition: second C and D})},\
\nonumber\\
{\rm and}\
L_n:&=&
e^{\xi_n}\Big(\sup_{n-1< t\leq n}\int_{(n-1)+}^t
e^{-\xi_{s-}}\ud\eta_s-\int_{(n-1)+}^n e^{-\xi_{s-}} \ud \eta_s\Big).
\end{eqnarray}
$(M_n,Q_n,L_n)_{n\in\bbn}$ is an iid sequence,
as follows by an easy extension of our proof of
Lemma~\ref{Appendix}. With these allocations $Z_n$ can be
written via (\ref{definition: second discrete series for Z})
in the form
\beam\label{Z}
Z_n=\sum_{i=1}^n\prod_{j=1}^{i-1}M_jQ_i=X_n-L_n\prod_{j=1}^nM_j.
\eeam
Nyrhinen proves the following result with equality in distribution:

\begin{prop}\label{pro}
Let $(M_n,Q_n,L_n)$ and $Z_n$ be as defined in
\eqref{definition: M allocation}, \eqref{definition: L allocation}
and \eqref{Z}.
Define $X_n$ as in \eqref{definition: Nyrhinens series}.
Then
\[\sup_{n-1<t\le n}Z_t=X_n\quad\mbox{and}\quad
\sup_{0\le t\le n} Z_t=\max_{m=1,\ldots,n}X_m.\]
\end{prop}

\bproof
For $n\in\bbn$ we have
\begin{eqnarray*}
\sup_{n-1< t\leq n}Z_t
&=& Z_{n-1}+\sup_{n-1< t\leq n}\int_{(n-1)+}^t e^{-\xi_{s-}}\ud\eta_s
\nonumber\\
&=& Z_{n-1}+
\int_{(n-1)+}^n e^{-\xi_{s-}}\ud\eta_s+e^{-\xi_n}L_n\\
& = &
X_n-\prod_{j=1}^n M_j L_n+e^{-\xi_n}L_n
\, = \,  X_n.
\end{eqnarray*}
This further implies that
$\sup_{0\le t\le n} Z_t=\max_{m=1,\ldots,n}X_m$.
\halmos
\eproof

Define the first passage time of $Z$ above $u>0$  by
$\tau_u^Z:=\inf\{t\ge 0:Z_t>u\}$.
Then Proposition \ref{pro} implies that for all $t>0$,
\[
P(\tau_u^Z\le t)=P(\tau_u^X\le t)
\quad {\rm and}\quad
P(\tau_u^Z<\infty)=P(\tau_u^X<\infty).
\]
So \eqref{4.15} and \eqref{equation: Nyrhinens finite time}
hold with $\tau_u^X$ replaced by $\tau_u^Z$,
when Hypothesis H is satisfied for the associated values of
$(M_n,Q_n,L_n).$ If, further, the distribution of $\ln M$ is spread
out, then \eqref{equation: Nyrhinens Goldie}
holds with $\tau_u^X$ replaced by $\tau_u^Z.$ This is the content of
Theorem 4 and Corollary 5 of  \cite{Nyrhinen01}.

\begin{rem} We make some comments on  Nyrhinen \cite{Nyrhinen01}.

(i) We begin with the discrete results.
Firstly, the sequence $X_n$ defined in (\ref{definition:
Nyrhinens series}) converges as $n\to\infty$ a.s. to a finite r.v. under Hypothesis H. To see this, note that if we
choose $L_n=L$ then $X_n$ is the inner iteration sequence $I_n(L)$
for the random equation
$\phi(t)=Mt+Q.$ Goldie and Maller \cite{goldie:maller:2000} prove that
$I_n(L)$ converges a.s. to a finite r.v. iff
$\prod_{j=1}^n M_j\rightarrow 0$ a.s. as $n\rightarrow\infty$ and
$I_{M,Q}<\infty,$ where $I_{M,Q}$ is an integral involving the
marginal distributions of $M$ and $Q.$ Since
these conditions have no
dependence on the distribution of $L$, it is clear that they
are precisely those under which $X_n$ converges a.s. for
iid $(M_n,Q_n,L_n).$ We now show that these conditions are in
fact satisfied under Hypothesis H, and thus the sequences $X_n$ and
$\sum_{i=1}^n\prod_{j=1}^{i-1}M_jQ_i$ converge a.s.,
and  to the same finite r.v..

Under Hypothesis H and  our assumption $P(M=0)=0$,
$E\ln M$ is well-defined and $E\ln M\in[-\infty,0)$.
Hence the random walk
$S_n:=\sum_{j=1}^n(-\ln M_j)=-\ln\prod_{j=1}^n M_j$ drifts to
$\infty$ a.s., and it follows that $\prod_{j=1}^n M_j\rightarrow 0$
a.s. as $n\rightarrow\infty$.
Since $\alpha_0^+>0$ there exists $s>0$ such that $E|Q|^s<\infty$,
thus $E\ln^+|Q|<\infty$. Hence Corollary~4.1 of
\cite{goldie:maller:2000} implies that the integral condition
$I_{M,Q}<\infty$ is satisfied and the sequence
$\sum_{i=1}^n\prod_{j=1}^{i-1}M_jQ_i$ converges a.s.

(ii) Nyrhinen transfers his discrete results into  continuous time, but
the corresponding results are difficult to apply in general.
The most problematic assumption is his condition
$\bar{y}=\infty$ (see \eqref{Nyybar}). In our notation, this is
equivalent to the condition $\psi(z)>0$ for all $z\ge 0.$ Theorem 1
of \cite{BankovskySly08} gives
necessary and sufficient
conditions on the \Levy measure of $(\xi,\eta)$ for this,
 which are amenable to verification in special cases,
as we showed in Section \ref{Examples}.
Verifying  Nyrhinen's condition $0<w^+<\alpha_0^+\le\infty$
requires finiteness of powers of $E|Z_1|$ and
$E[\sup_{0<t\le 1}|Z_t|]$.
These conditions would be more conveniently
stated in terms of the characteristic
triplet of $(\xi,\eta)$ or (at least) the marginal distributions of
$\xi$ and $\eta.$
In the special case that $\xi$ and $\eta$ are
independent \Levy processes, Theorem 3.2 of Paulsen
\cite{Paulsen02} does exactly that. However, problems remain.
In  \cite{Paulsen02}, the condition $\bar y=\infty$ is assumed to
be true whenever  $\xi$ and $\eta$ are independent and $\eta$ is not
a subordinator. However, this claim is false$^\dagger$. \footnote{$^\dagger$To see this,
let $(\xi,\eta)_t:=(t+N_t,-t)$ where
$N$ is a Poisson process with jump times $0<\tau_0<\tau_1<\cdots$. This
example trivially satisfies all the conditions in Paulsen's Theorem~3.2. However, using Ito's formula for semi-martingales and some
simple manipulation we obtain
$Z_t=-1+(e-1)\sum_{i=1}^{N_t}e^{-\tau_i-i}+e^{-t-N_t}$, and hence
$\inf_{t>0}Z_t\ge-1$ a.s.}
(It does hold if extra conditions are
imposed, in line with Remark 2(3) of \cite{BankovskySly08}.)
Finally, it would be desirable to remove the finite mean assumption
for $\xi$ in \cite{Paulsen02}
and replace the  moment conditions in \cite{Paulsen02},
which are sufficient
for convergence of $Z_t$, with the precise necessary and sufficient
conditions given in Goldie and Maller~\cite{goldie:maller:2000}.
Our Theorem \ref{theorem: my Nyrhinen result} addresses all of the
above concerns in the most general setting.
\end{rem}

\section{Proof of Theorem \ref{theorem: my Nyrhinen result}}
The proof requires the
following lemma, which was stated but not proved in \cite{Bankovsky09}.

\begin{lem}\label{sup proposition} Suppose there exist
$r>0$ and $p,q>1$ with $1/p+1/q=1$ such that $E e^{-\max \{1,r
\}p\xi_1} <\infty$ and $E|\eta_1| ^{\max\{1,r
\}q}<\infty.$
Then
\begin{equation}\label{equation: sup proposition}
E\Big[\sup_{0\le t\le 1}\left|Z_t\right|^{\max\{1,r \}}\Big]
=
E\Big[\sup_{0\le t\le 1}
\Big|\int_0^t e^{-\xi_{s-}}\ud\eta_s\Big|^{\max\{1,r \}}\Big]
<\infty.
\end{equation}
\end{lem}

\bproof
For ease of notation let $k:=\max\{1,r\}.$ Assume there exists $r>0$
and $p,q>1$ with $1/p+1/q=1$ such that
$Ee^{-kp\xi_1}<\infty$ and $E|\eta_1|^{kq}<\infty.$
We prove the lemma first for the case
in which $E\eta_1=0.$ Since $\eta$ is a \Levy process this implies that $\eta$ is a \cadlag martingale. Since $\xi$
is \cadlag $e^{-\xi}$ is a locally bounded process and hence $Z$ is
a local martingale for $\mathbb{F}$ by the construction of the stochastic integral (see e.g. Protter~\cite{protter}). Since additionally $Z_0=0,$ the Burkholder-Davis-Gundy inequalities
ensure that for our choices of $p,q$ and $k$ there
exists $b>0$ such that
\begin{eqnarray*}E\Big[\sup_{0\le t\le 1}\Big|\int_0^t
e^{-\xi_{s-}}\ud\eta_s\Big|^k\Big]
&\le& bE\Big[\Big[\int_0^z e^{-\xi_{s-}}\ud\eta_s,\int_0^z e^{-\xi_{s-}}\ud\eta_s\Big]_{z=1}^{z=k/2}\Big]\\
= \, bE\Big[\Big(\int_0^1 e^{-2\xi_{s-}}\ud
[\eta,\eta]_s\Big)^{k/2}\Big]
&\le& bE\Big[\Big(\int_0^1 \sup_{0\le t\le 1}e^{-2\xi_t}\ud
[\eta,\eta]_s\Big)^{k/2} \Big],
 \end{eqnarray*}
where in the second inequality recall that
$[\eta,\eta]_s$ is increasing.
(The notation $[\cdot,\cdot]$ denotes the quadratic variation process.)
The last expression equals
\[
bE\Big[\sup_{0\le t\le 1}e^{-k\xi_{t}}[\eta,\eta]_1^{k/2}\Big]
\le b\Big(E\Big[\sup_{0\le t\le 1}e^{-pk\xi_{t}}\Big]
\Big)^{1/p}\Big(E\big[[\eta,\eta]_1^{qk/2}\big]\Big)^{1/q},
\]
where the inequality follows for
our choices of $p$ and $q$ by H\"{o}lder's inequality.
Since $k\ge1$, $q>1$, the Burkholder-Davis-Gundy inequalities give
the existence of $c>0$ such that (using Doob's inequality for the second inequality)
\begin{eqnarray*}
E\left[[\eta,\eta]_1^{qk/2}\right] \leq \frac{1}{c}E\Big[\sup_{0\le
t\le 1}|\eta_t|^{qk}\Big]
\le\frac{8}{c}E\left[|\eta_1|^{qk}\right]<\infty.
\end{eqnarray*}

Thus it suffices to prove
$E\left[\sup_{0\le t\le 1}e^{-pk\xi_{t}}\right]<\infty$.
Now $Y_t:=e^{-pk\xi_t}/c^t$,  where
$c:=Ee^{-pk\xi_1} \in(0,\infty)$
is a non-negative martingale, and
it follows by Doob's maximal inequality that
\begin{equation*}
E\Big[\sup_{0\le t\le1}e^{-pk\xi_{t}}\Big]\le\max\{1,c\}
E\Big[\sup_{0\le t\le 1}
\frac{e^{-pk\xi_t}}{c^t} \Big]
\le
\max\Big\{\frac1{c},1\Big\}
\Big(\frac{pk}{pk-1}\Big)^{pk}
Ee^{-pk\xi_1}<\infty.
\end{equation*}
Hence the lemma is proved for the case in which $E(\eta_1)=0.$
In general, write
\begin{eqnarray*}
&&E\Big[\sup_{0 \le t\le
1}\Big|\int_0^t e^{-\xi_{s-}}\ud \eta_s\Big|^k\Big]
= E\Big[\sup_{0 \le t\le 1}\Big|\int_0^t e^{-\xi_{s-}}\ud(
\eta_s-sE\eta_1+sE\eta_1)\Big|^k
\Big]\\
&&\le E\Big[\Big(\sup_{0 \le t\le 1}\Big|\int_0^t
e^{-\xi_{s-}}\ud( \eta_s-sE\eta_1)\Big|+
|E\eta_1|\sup_{0 \le t\le 1}\Big|\int_0^t
e^{-\xi_{s-}}\ud s\Big| \Big)^k\Big],
\end{eqnarray*}
in which the first term on the right-hand side is finite by the first
part of the proof.
An application of Minkowski's inequality to the second term on the right-hand side completes the proof.
\halmos
\eproof

\begin{rem}\rm
If $\xi$ and $\eta$ are independent, then H\"{o}lder's
inequality is not required in the proof of Lemma \ref{sup
proposition}, and a simpler independence argument shows that
(\ref{equation: sup proposition}) holds if $Ee^{-\max
\{1,r \}\xi_1}<\infty$ and $E|\eta_1| ^{\max\{1,r
\}}<\infty$ for some $r>0.$ We can put further restrictions
on $\xi$ and $\eta,$ such as in the example in Section 3 of Nyrhinen
\cite{Nyrhinen01}, which assumes $\xi$ is continuous and $\eta$ is
compound Poisson plus drift, which render the use of the
Burkholder-Davis-Gundy inequalities unnecessary and further simplify
the conditions. For general \Levy $(\xi,\eta)$ the above inequality
is the sharpest we have found.
\end{rem}

{\em Proof of Theorem \ref{theorem: my Nyrhinen result}: }
We aim to use  Proposition  \ref{theorem: Nyrhinens main theorem}
for passage below rather than above.
We can do this by replacing $\eta$ by $-\eta$. Note that for $z>0$,
\[
T_z=\inf\{t>0:Z_t<-z\} =\inf\{t>0:-Z_t>z\}
=\inf\{t>0:\whZ_t>z\},
\]
where we denote $Z_t$, when $\eta$ is replaced by $-\eta$, by $\whZ_t$
and similarly for the other quantities.
Thus $\whZ_t=-Z_t$, and it is easily checked that,
with $(M_n,Q_n)$ as in \eqref{definition: M allocation},
$(\whM_n,\whQ_n)= (M_n,-Q_n)$, and,
with $L_n$ as in \eqref{definition: M allocation},
$\whL_n=-\barL_n$, where
\begin{equation}\label{Ltilde}
\barL_n:=-e^{\xi_n}\Big(
\int_{(n-1)+}^n e^{-\xi_{s-}} \ud \eta_s-
\inf_{n-1< t\leq n}\int_{(n-1)+}^te^{-\xi_{s-}}\ud\eta_s
\Big).
\end{equation}
{}From \eqref{definition: Nyrhinens series} we get
$\whX_n(\whL_n)=-X_n(\barL_n)$.
Then  Proposition  \ref{theorem: Nyrhinens main theorem}
ensures that
(\ref{equation: Nyrhinen result 1}) and
(\ref{equation: Nyrhinen result 2})
hold, if we can
prove that the relevant conditions are satisfied for $(\whM,\whQ,\whL)$;
i.e., we must show that Hypothesis H holds for the hat variables.

The corresponding
$\wh {\overline y}$ (see \eqref{Nyybar}) is
\begin{eqnarray*}
\sup\Big\{y\in\mathbb{R}:
P\Big(\sup_{n\in\mathbb{N}}\whX_n(\whL_n)>y\Big)>0\Big\}
&=&
\inf\Big\{z\in\mathbb{R}:
P\Big(\inf_{t>0}Z_t<-z\Big)>0\Big\},
\end{eqnarray*}
so  $\wh {\overline y}=\infty$
if and only if $\psi(z)>0$ for all $z\ge 0,$ which we have assumed.

We need a $w^+>0$ such that $E\whM^{w^+}=1$, and this is the case
with $w^+=w$ under
\eqref{cramer} since $\whM=M=e^{-\xi_1}$.
Also, $\whc_M(\alpha)= \ln E\whM^\alpha=c(\alpha)$, so that
$\alpha_0^+$ in \eqref{t0+} here equals
$ \alpha_0$ as defined in \eqref{t0}.
Note that the extra term
$E(\whM\whL^+)^\alpha=EM\barL^-)^\alpha$
required in \eqref{t0} is superfluous here, since
$E(M\barL^-)^\alpha= E\barZ_1^\alpha$, and this is finite for
$\alpha\ge 0$ if and only if $E|Z_1|^\alpha<\infty$.

Under the  moment conditions of
Theorem \ref{theorem: my Nyrhinen result},
the conditions of Lemma \ref{sup proposition} hold with $r=w^++\veps$,
so  $E|Z_1|^\alpha<\infty$ for $\alpha=\max\{1,w+\veps\}$,
and hence $\alpha_0^+\ge w^++\veps>w^+$.
Thus indeed Hypothesis H is fulfilled in the present situation and
Proposition  \ref{theorem: Nyrhinens main theorem}
applies to give
(\ref{equation: Nyrhinen result 1}) and
(\ref{equation: Nyrhinen result 2}).
Also  $\alpha_0^+\ge w^++\veps>w^+$  implies
$c'(\alpha_0-)>c'(w)=\mu^*=-E\xi_1e^{-w\xi_1}$,
and this is finite since $Ee^{-(w+\veps)\xi_1}$ is.
So $0\le  \alpha_0<1/\mu^*<\infty$.

Suppose, further, that $\xi_1$ is spread out. Then the dual
version of (\ref{equation: Goldie result}) follows
from Nyrhinen's comments in \cite{Nyrhinen01}, which we expressed as Proposition \ref{theorem: Nyrhinens main theorem}.
 \halmos

\end{document}